\documentclass[12pt,a4paper]{amsart}
\usepackage{amsmath,amsxtra,amsfonts,amssymb,
latexsym,amsthm,amscd,verbatim,amstext,enumerate,graphicx}
\usepackage[mathscr]{eucal}
\makeatletter
\@namedef{subjclassname@2010}{%
  \textup{2010} Mathematics Subject Classification}
\makeatother
\newcommand{\me}{\mathrm{e}}

\newcommand{\dif}{\mathrm{d}}

\newtheorem{Th}{Theorem}[section]

\newtheorem{Lem}[Th]{Lemma}
\newtheorem{Md}[Th]{Proposition}

\theoremstyle{definition}

\newtheorem{Not}[Th]{Remark}

\numberwithin{equation}{section}

\begin{document}

%%%%% To ease editing, for IMPAN journals add:

\baselineskip=17pt

%%%%%%%%%%%

%% In the running head, replace first names by initials 
%% and give an abbreviation of the title.

\title[]{On the initial value problem for the Navier-Stokes equations with the initial datum in the Sobolev spaces}

\author[Dao Quang Khai]{D. Q. Khai}
\address{Institute of Mathematics\\ Vietnam Academy of Science and Technology\\
18 Hoang Quoc Viet, 10307  Cau Giay, Hanoi, Vietnam.}
\email{khaitoantin@gmail.com}

\author[Vu Thi Thuy Duong]{V. T. T. Duong}
\address{Faculty of Basic Sciences \\ Quang Ninh University of Industry\\
Yen Tho, Dong Trieu, Quang Ninh, Vietnam.}
\email{vuthuyduong309@gmail.com}

\date{}

\begin{abstract}
In this paper, we study local well-posedness for the Navier-Stokes equations with arbitrary initial data in homogeneous Sobolev spaces $\dot{H}^s_p(\mathbb{R}^d)$ for $d \geq 2, p > \frac{d}{2},\ {\rm and}\ \frac{d}{p} - 1 \leq s  < \frac{d}{2p}$. The obtained result improves the known ones for $p > d$ and $s = 0$ (see \cite{M. Cannone 1995, M. Cannone. Y. Meyer 1995}). In the case of critical indexes $s=\frac{d}{p}-1$, we prove global well-posedness for Navier-Stokes equations when the norm of the initial value is small enough. This result is a generalization of the ones in  \cite{M. Cannone 1999} and \cite{P. G. Lemarie-Rieusset 2002}  in which $(p = d, s = 0)$ and $(p > d, s = \frac{d}{p} - 1)$, respectively.
\end{abstract}



\maketitle

\section{Introduction}

This paper studies the Cauchy problem of the incompressible Navier-Stokes equations (NSE)
in the whole space $\mathbb R^d$ for $d\geq 2$,
\begin{align} 
\left\{\begin{array}{ll} \partial _tu  = 
\Delta u - \nabla .(u \otimes u) - \nabla p , & \\ 
\nabla .u = 0, & \\
u(0, x) = u_0, 
\end{array}\right . \notag
\end{align}
which is a condensed writing for
\begin{align} 
\left\{\begin{array}{ll} 1 \leq k \leq d, \ \  \partial _tu_k  
= \Delta u_k - \sum_{l =1}^{d}\partial_l(u_lu_k) - \partial_kp , & \\ 
\sum_{l =1}^{d}\partial_lu_l = 0, & \\
1 \leq k \leq d, \ \ u_k(0, x) = u_{0k} .
\end{array}\right . \notag
\end{align}
The unknown quantities are the velocity 
$u(t, x)=(u_1(t, x),\dots,u_d(t, x))$ of the fluid 
element at time $t$ and position $x$ and the pressure $p(t, x)$.\\
There is an extensive literature on the existence of strong solutions of the Cauchy problem for NSE. The global well-posedness of strong  solutions for small initial data in the critical Sobolev space  $\dot{H}^{\frac{1}{2}}(\mathbb R^3)$ is due to Fujita and Kato \cite{T. Kato 1962}, also in \cite{J. M. Chemin 1992}, Chemin has proved the case of  $H^s(\mathbb R^3), (s >1/2)$. In \cite{T. Kato 1984}, Kato has proved the case of the Lebesgue space $L^3(\mathbb R^3)$. In \cite{H. Koch 2001}, Koch and Tataru have proved the case of the space $BMO^{-1}$ (see also \cite{J. Y. Chemin 2009}). In \cite{H. Koch 2001}, H. Koch has proved the case of the space $\dot{B}^{\frac{d}{p}-1}_{p, \infty}(\mathbb{R}^d)_{(p<+\infty)}$, see \cite{H. Koch 2001} and the recent ill-posedness result  \cite{J. Bourgain 2008} for $\dot{B}^{- 1}_{\infty, \infty}(\mathbb{R}^d)$. Results on the existence of mild solutions with value in  $L^p(\mathbb{R}^d), (p > d)$ were established in  the papers of  Fabes, Jones and Rivi\`{e}re \cite{E. Fabes 1972a} 
and of Giga \cite{Y. Giga 1986a}. Concerning the initial 
datum in the space $L^\infty$, the existence of a mild solution was 
obtained by Cannone and Meyer in \cite{M. Cannone 1995,M. Cannone. Y. Meyer 1995}. Moreover, in \cite{M. Cannone 1995,M. Cannone. Y. Meyer 1995}, they also obtained theorems 
on the existence of mild solutions with value in the Morrey-Campanato space 
$M^p_2(\mathbb{R}^d), (p > d)$ and the 
Sobolev space $H^s_p(\mathbb{R}^d), (p < d, \frac{1}{p} 
- \frac{s}{d} < \frac{1}{d})$. NSE in the Morrey-Campanato space were also treated by Kato  \cite{T. Kato 1992} and Taylor \cite{M.E. Taylor 1992}. Recently, the authors of this article have considered NSE in Sobolev spaces, Sobolev-Lorentz spaces, mixed-norm Sobolev-Lorentz spaces, and Sobolev-Fourier-Lorentz spaces, see \cite{N. M. Tri: Tri2014???,N. M. Tri: Tri2014????,N. M. Tri: Tri2016??}, \cite{N. M. Tri: Tri2015??}, \cite{N. M. Tri: Tri2014a}, and \cite{N. M. Tri: Tri2015?} respectively. In \cite{N. M. Tri: Tri2016?}, we prove some results on the existence and  space-time decay rates of global strong solutions of the Cauchy problem for NSE in weighed $L^\infty(\mathbb R^d,|x|^\beta{\rm dx})$ spaces. In this paper, we construct mild solutions in the spaces 
$L^\infty([0, T];\dot{H}^s_p(\mathbb{R}^d))$ 
to the Cauchy problem for NSE when  the initial datum belongs to 
the Sobolev spaces $\dot{H}^s_p(\mathbb{R}^d)$, 
with $d \geq 2, p > \frac{d}{2},\ {\rm and}\ \frac{d}{p} 
- 1 \leq s  < \frac{d}{2p}$, we obtain the existence of mild solutions with  
arbitrary initial value when $T$ is small enough  and existence of mild 
solutions for any $T < +\infty$ when the norm of the initial value in 
the Triebel-Lizorkin spaces $\dot{F}^{s- d(\frac{1}{p} 
- \frac{1}{\tilde q}), \infty}_{\tilde q}$, $(\tilde q > 
{\rm max}\{p, q\},\ \text{where}\ \frac{1}{q} = \frac{1}{p} - \frac{s}{d})$ 
is small enough. In the case $p > d$ and $s = 0$, this result is stronger than that of Cannone and Meyer \cite{M. Cannone 1995, M. Cannone. Y. Meyer 1995} under a weaker condition on the initial data. In the case of critical indexes $(p > \frac{d}{2}, 
s = \frac{d}{p} - 1)$, we obtaine global mild solutions when the norm of the initial value in the Triebel-Lizorkin spaces $\dot{F}^{\frac{d}{\tilde q} - 1, \infty}_{\tilde q}(\mathbb{R}^d), (\tilde q > {\rm max}\{d, p\})$ is small enough. This result in one hand if $p = d$ and $s = 0$ is stronger than that of Cannone and Planchon \cite{M. Cannone 1999} but under a weaker condition on the initial data  and in the other hand if $p > d$ and $s = \frac{d}{p} - 1$ is stronger than that of Lemarie-Rieusset but under a weaker condition on the initial data (Proposition 20.2, \cite{P. G. Lemarie-Rieusset 2002}, p. 201). The content of this paper is as follows: in Section 2, we state our main theorem after introducing some notations. In Section 3, we first establish some estimates concerning the heat semigroup with differential. We also recall some auxiliary lemmas and several estimates in the homogeneous Sobolev spaces and Triebel spaces. Finally, in Section 4, we will give the proof of the main theorem.
\section{Statement of the results}
Now, for $T > 0$, we say that $u$ is a mild solution of NSE on $[0, T]$ 
corresponding to a divergence-free initial datum $u_0$ 
when $u$ solves the integral equation
$$
u = e^{t\Delta}u_0 - \int_{0}^{t} e^{(t-\tau) \Delta} \mathbb{P} 
\nabla  .\big(u(\tau,.)\otimes u(\tau,.)\big) \dif\tau.
$$
Above we have used the following notation: 
for a tensor $F = (F_{ij})$ 
we define the vector $\nabla.F$ by 
$(\nabla.F)_i = \sum_{j = 1}^d\partial_jF_{ij}$ 
and for two vectors $u$ and $v$, we define their tensor 
product $(u \otimes v)_{ij} = u_iv_j$. 
The operator $\mathbb{P}$ is the Helmholtz-Leray 
projection onto the divergence-free fields 
$$
(\mathbb{P}f)_j =  f_j + \sum_{1 \leq k \leq d} R_jR_kf_k, 
$$
where $R_j$ is the Riesz transforms defined as 
$$
R_j = \frac{\partial_j}{\sqrt{- \Delta}}\ \ {\rm i.e.} \ \  
\widehat{R_jg}(\xi) = \frac{i\xi_j}{|\xi|}\hat{g}(\xi).
$$
The heat kernel $e^{t\Delta}$ is defined as 
$$
e^{t\Delta}u(x) = ((4\pi t)^{-d/2}e^{-|.|^2/4t}*u)(x).
$$
For a space of functions defined on $\mathbb R^d$, 
say $E(\mathbb R^d)$, we will abbreviate it as $E$.
We denote by $L^q: = L^q(\mathbb R^d)$ the usual 
Lebesgue space for $q \in [1, \infty]$ with the norm $\|.\|_q$, 
and we do not distinguish between the vector-valued and scalar-valued spaces of functions.
 We define the Sobolev space by $\dot{H}^s_q: = \dot{\Lambda}^{-s}L^q$ 
equipped with the norm 
$\big\|f\big\|_{\dot{H}^s_q}: =  \|\dot{\Lambda}^sf\|_q$. 
Here $\dot{\Lambda}^s:= \mathcal F^{-1}|\xi|^s\mathcal F$, 
where $\mathcal F$ and $\mathcal F^{-1}$ are the Fourier transform and its inverse, respectively. $\dot{\Lambda} = \sqrt{-\Delta}$ is the homogeneous Calderon pseudo-differential operator. For vector-valued $f=(f_1,...,f_M)$, we define $\|f\|_X = \big(\sum_{m=1}^{m=M}\|f_m\|_X^2\big)^{\frac{1}{2}}$. Throughout the paper, we sometimes use the notation $A \lesssim B$ as an equivalent to $A \leq CB$ with a uniform constant $C$. The notation $A \simeq B$ means that $A \lesssim B$ and $B \lesssim A$. Now we can state our results
\begin{Th}\label{th1} Let $s$ and $p$ be such that
$$
p > \frac{d}{2}\ \text{and}\ \frac{d}{p} 
- 1 \leq s  < \frac{d}{2p}.
$$
Set
$$
\frac{1}{q} = \frac{1}{p} - \frac{s}{d}.
$$
{\rm (a)} For all $\tilde q > {\rm max}\{p, q\}$, 
there exists a positive constant $\delta_{q,\tilde q,d}$ 
such that for all $T > 0$ and for all 
 $u_0 \in \dot{H}^s_p(\mathbb{R}^d)\ 
with\ {\rm div}(u_0) = 0$ satisfying
\begin{equation}\label{eq1}
T^{\frac{1}{2}(1+s-\frac{d}{p})}
\Big\|\underset{0 < t < T}{\rm sup}t^{\frac{d}{2}(\frac{1}{p}
- \frac{s}{d} - \frac{1}{\tilde q})}
\big|\me^{t\Delta}u_0\big|\Big\|_{L^{\tilde q}} \leq \delta_{q,\tilde q,d},
\end{equation}
NSE has a unique mild solution $u \in L^\infty([0, T]; \dot{H}^s_p)$ and the following inequality holds
$$
\Big\|\underset{0 < t < T}{\rm sup}t^{\frac{d}{2}(\frac{1}{q}-\frac{1}{r})}\big|u(t,x)\big|\Big\|_r < +\infty,\ for\ all\ r > {\rm max}\{p, q\}.
$$
In particular, the condition \eqref{eq1} holds for arbitrary  $u_0 \in \dot{H}^s_p(\mathbb{R}^d)$ when  $T(u_0)$ is small enough.\\
{\rm (b)} If $s = \frac{d}{p} - 1$ then for all 
$\tilde q  > {\rm max}\{p,d\}$ there exists a constant 
$\sigma_{\tilde q, d}>0$ such that if  
$\big\|u_0\big\|_{\dot{F}^{\frac{d}{\tilde q} - 1, \infty}_{\tilde q}} 
\leq \sigma_{\tilde q,d}$ and $T = +\infty$ 
then the condition \eqref{eq1} holds.
\end{Th}

In the case of critical indexes $(s = \frac{d}{p} - 1, 
p > \frac{d}{2})$, we get the following consequence.
\begin{Md}\label{pro1} Let $ p > \frac{d}{2}$. 
Then for any  $\tilde q > {\rm max}\{p, d\}$, 
there exists a positive constant $\delta_{\tilde q,d}$ 
such that for all $T > 0$ and for all 
 $u_0 \in \dot{H}^{\frac{d}{p} - 1}_p(\mathbb{R}^d)$ 
with ${\rm div}(u_0) = 0$ satisfying
\begin{equation}\label{eq2}
\Big\|\underset{0 < t < T}{\rm sup}
t^{\frac{1}{2}(1 - \frac{d}{\tilde q})}
\big|\me^{t\Delta}u_0\big|\Big\|_{L^{\tilde q}} 
\leq \delta_{\tilde q,d},
\end{equation}
NSE has a unique mild solution $u \in L^\infty([0, T]; \dot{H}^{\frac{d}{p} - 1}_p)$ and the following inequality holds
$$
\Big\|\underset{0 < t < T}{\rm sup}t^{\frac{d}{2}(\frac{1}{d}-\frac{1}{r})}\big|u(t,x)\big|\Big\|_r < +\infty,\ for\ all\ r > {\rm max}\{p, d\}.
$$
Denoting $w = u - \me^{t\Delta}u_0$ then $w$ satisfies the following inequality
$$
\Big\|\underset{0 < t < T}{\rm sup}
\big|\dot{\Lambda}^{\frac{d}{\tilde p}-1}w(t,x)\big|\Big\|_{L^{\tilde p}} < \infty, \ for\ all\ \tilde p >  \frac{1}{2}{\rm max}\{p, d\}.
$$
Moreover, if $p \geq d$ then 
$$
\underset{0 < t < T}{\rm sup}
t^{\frac{1}{2}({1-\frac{d}{p}})}\big|u(t, .)\big| \in L^p.
$$
In particular, the condition \eqref{eq2} 
holds for arbitrary  $u_0 \in \dot{H}^{\frac{d}{p} - 1}_p(\mathbb{R}^d)$ 
when  $T(u_0)$ is small enough, and there exists a positive 
constant $\sigma_{\tilde q,d}$ such that if  
\begin{equation}\label{eq3}
\big\|u_0\big\|_{\dot{F}^{\frac{d}{\tilde q} - 1, \infty}_{\tilde q}} 
\leq \sigma_{\tilde q,d}\ and\ T = +\infty
\end{equation} then the condition \eqref{eq2} holds. 
\end{Md}
\begin{Not} Proposition \ref{pro1} is the theorem of Canone and Planchon \cite{M. Cannone 1999} if $p = d$ and the condition \eqref{eq3} is replaced by 
the condition 
\begin{equation}\label{eq4}
\big\|u_0\big\|_{\dot{F}^{\frac{d}{\tilde q} - 1, \infty}_{\tilde q}}
 \leq \sigma_{\tilde q,d} \ {\rm and}\ T = +\infty,\ {\rm where }\ d <\tilde q <2d.
\end{equation}
Note that in the case $p=d$, the condition \eqref{eq3} is weaker than 
the condition \eqref{eq4} because of the 
following elementary imbedding maps
$$
\dot{F}^{\frac{d}{\tilde q}-1, \infty}_{\tilde q}
(\mathbb{R}^d)_{(d <\tilde q <2d)}
\hookrightarrow  \dot{F}^{-\frac{1}{2}, 
\infty}_{2d}(\mathbb{R}^d) \hookrightarrow  
\dot{F}^{\frac{d}{\tilde q}-1, 
\infty}_{\tilde q}(\mathbb{R}^d)_{(\tilde q > 2d)}.
$$
\end{Not}
\begin{Not} The statement about the global existence in Proposition \ref{pro1} is the Lemarie-Rieusset statement 
(Proposition 20.2, \cite{P. G. Lemarie-Rieusset 2002}, p. 201) if  $p > d$ and the condition \eqref{eq3} is replaced by 
the condition 
\begin{equation}\label{eq5}
\big\|u_0\big\|_{\dot{H}^{\frac{d}{p} - 1}_p} < \delta_{d,p}.
\end{equation}
Note that the condition \eqref{eq3} is weaker than the 
condition \eqref{eq5} because of the following 
elementary imbedding maps
$$
\dot{H}^{\frac{d}{p}-1}_p(\mathbb{R}^d) 
\hookrightarrow  \dot{F}^{\frac{d}{p}-1, \infty}_p(\mathbb{R}^d) 
\hookrightarrow  \dot{F}^{\frac{d}{\tilde q}-1, \infty}_{\tilde q}
(\mathbb{R}^d), (\tilde q > p).
$$
Lemarie-Rieusset proved the above statement by using 
Hardy-Littlewood maximal functions theory 
(as developped for $L^d$ by Cander\'{o}n \cite{C. Calderon 1993}  
and \linebreak Cannone \cite{M. Cannone 1995}).
\end{Not}
In the case of supercritical indexes $p > \frac{d}{2}\ \text{and}\ \frac{d}{p} - 1 < s  < \frac{d}{2p}$, we get the following consequence.
\begin{Md}\label{pro2}
Let $p > \frac{d}{2}\ \text{and}\ \frac{d}{p} - 1 < 
s  < \frac{d}{2p}$. Then for all $\tilde q > {\rm max}\{p, q\}$, 
where
$$
\frac{1}{q} = \frac{1}{p} - \frac{s}{d},
$$
there exists a positive constant $\delta_{q,\tilde q,d}$ 
such that for all $T > 0$ and for all \linebreak
 $u_0 \in \dot{H}^s_p(\mathbb{R}^d)\ with\ {\rm div}(u_0) = 0$ satisfying
\begin{equation}\label{eq6}
T^{\frac{1}{2}(1+s-\frac{d}{p})}
\big\|u_0\big\|_{\dot{F}^{s- (\frac{d}{p} - \frac{d}{\tilde q}),
 \infty}_{\tilde q}} \leq \delta_{q,\tilde q,d},
\end{equation}
NSE has a unique mild solution $u \in L^\infty([0, T]; \dot{H}^s_p)$ and the following inequality holds
$$
\Big\|\underset{0 < t < T}{\rm sup}t^{\frac{d}{2}(\frac{1}{q}-\frac{1}{r})}\big|u(t,x)\big|\Big\|_r < +\infty,\ for\ all\ r > {\rm max}\{p, q\}.
$$
\end{Md}
\begin{Not}
Proposition \ref{pro2} is the theorem of Canone and Meyer \cite{M. Cannone 1995,M. Cannone. Y. Meyer 1995} if $s=0$, $p > d$, and the condition \eqref{eq6} is replaced by the condition
\begin{equation*}
T^{\frac{1}{2}(1-\frac{d}{p})}\big\|u_0\big\|_{L^p} \leq \delta_{p,d}.
\end{equation*}
Note that in the case $s=0$ and $p > d$, the condition \eqref{eq6} is weaker than the above condition because of the following elementary imbedding maps
$$
L^p(\mathbb{R}^d) \hookrightarrow  \dot{F}^{-(
\frac{d}{p}-\frac{d}{\tilde q}), \infty}_{\tilde q}(\mathbb{R}^d), (\tilde q > p \geq d).
$$
\end{Not}
\section{Tools from harmonic analysis}
In this section we prepare some auxiliary lemmas.\\
The main property we use throughout this paper is that the operator $e^{t\Delta}\mathbb{P}\nabla$ is a matrix of convolution operators with bounded integrable kernels.
\begin {Lem} \label {lem1} Let $s>-1$. Then the kernel function of  $\dot{\Lambda}^se^{t\Delta} \mathbb{P}\nabla$ is the function
$$
K_t(x)= \frac{1}{t^{\frac{d+s+1}{2}}}K\big(\frac{x}{\sqrt t}\big),
$$
where the function $K$ is the kernel function of  $\dot{\Lambda}^se^{\Delta} \mathbb{P}\nabla$ which satisfies the following inequality
$$
|K(x)| \lesssim \frac{1}{1+|x|^{d+s+1}}.
$$
\end{Lem}
\begin{proof}See Proposition 11.1 in (\cite{P. G. Lemarie-Rieusset 2002}, p. 107).  \end{proof}
\begin{Lem}\label{lem2}
The kernel function $K_t(x)$ of  $\dot{\Lambda}^se^{t\Delta} \mathbb{P}\nabla$ satisfies the following inequality
$$
|K_t(x)| \lesssim \frac{1}{t^{\gamma_2}|x|^{\gamma_1}}, \  for\  \gamma_1 >0, \gamma_2 > 0,\ and\ \gamma_1 + 2\gamma_2 = d+s+1.
$$
\end{Lem}
\begin{proof}This is deduced by applying Lemma \ref{lem1} and the Young inequality
\begin{gather*}
|K_t(x)| = \Big|\frac{1}{t^{\frac{d+s+1}{2}}}
K\Big(\frac{x}{\sqrt t}\Big)\Big| 
\lesssim \frac{1}{t^{\frac{d+s+1}{2}}}
\frac{1}{1 + (\frac{|x|}{\sqrt t})^{d+s+1}} \notag\\
= \frac{1}{t^{\frac{d+s+1}{2}}+ |x|^{d+s+1}} 
\lesssim \frac{1}{{t}^{\gamma_2}|x|^{\gamma_1}}.
\end{gather*}
\end{proof}
\begin{Lem}\label{lem3} {\rm(Sobolev inequalities)}.\\
If \ $s_1 > s_2, \ 1 < q_1, \ q_2 < \infty$, and $s_1 - \frac{d}{q_1} = s_2 - \frac{d}{q_2}$, then we have the following embedding mapping
$$
\dot{H}_{q_1}^{s_1} \hookrightarrow \dot{H}_{q_2}^{s_2}.
$$
\end{Lem}
In this paper we use the definition of the homogeneous Triebel space  $\dot F^{s,p}_q$ in \cite{G. Bourdaud 1993,G. Bourdaud 1988,M. Frazier 1991,J. Peetre 1976}. The following lemma will provide a different characterization 
of Triebel spaces $\dot F^{s,p}_q$ in terms of the heat semigroup 
and will be one of the staple ingredients of the proof of Theorem \ref{th1}.
\begin{Lem}\label{lem4} \ \\
Let $1 \leq p, q \leq \infty$ and $ s < 0$. Then the two quantities
\begin{gather*}
\Big\|\Big(\int_0^\infty\big(t^{-\frac{s}{2}}\big|\me^{t\Delta}f\big|\big)^p
\frac{{\rm d}t}{t}\Big)^{\frac{1}{p}}\Big\|_{L^q}\ and 
\ \big\|f\big\|_{\dot{F}_{q}^{s, p}} \ are \ equivalent.
\end{gather*}
\end{Lem}
\begin{proof}See \cite{M. Cannone 1999}. \end{proof}
\begin{Lem}\label{lem5}{\rm (Convolution of  the Lorentz spaces).}\\
Let $1 < p <\infty $, $1 \leq q \leq \infty,  1/p' + 1/p = 1$, and $1/q' + 1/q = 1$. Then convolution is a bounded bilinear operator:\\
{\rm (a)} from $L^{p, q} \times L^1$ to $L^{p, q}$,\\
{\rm (b)} from $L^{p, q} \times L^{p', q'}$ to $ L^\infty$,\\
{\rm (c)} from $L^{p, q} \times L^{p_1, q_1}$ to $L^{p_2, q_2}$, for $ 1 < p, p_1, p_2 <\infty, 1 \leq q, q_1, q_2 \leq \infty,
 1/p_2 + 1 = 1/p + 1/p_1$, and  $1/q_2 = 1/q + 1/q_1$.
\end{Lem}
\begin{proof}See Proposition 2.4 (c) in (\cite{P. G. Lemarie-Rieusset 2002}, p. 20). \end{proof}
\begin{Lem}\label{lem6} Let $\theta< 1$ and $\gamma < 1$ then
$$
\int^t_0(t-\tau)^{-\gamma} \tau^{-\theta}{\rm d}\tau = C t^{1-\gamma - \theta},\  where\  
C = \int^1_0(1-\tau)^{-\gamma} \tau^{-\theta}{\rm d}\tau < \infty.
$$
\end{Lem}
The proof of this lemma is elementary and may be omitted.\qed \\
Let us recall following result on solutions of a quadratic
equation in Banach spaces (Theorem 22.4 in (\cite{P. G. Lemarie-Rieusset 2002}, p. 227)).
\begin{Th}\label{th2} Let $E$ be a Banach space, and $B: E \times E \rightarrow  E$ 
be a continuous bilinear map such that there exists $\eta > 0$ so that
$$
\|B(x, y)\| \leq \eta \|x\| \|y\|,
$$
for all x and y in $E$. Then for any fixed $y \in E$ 
such that $\|y\| \leq \frac{1}{4\eta}$, the equation $x = y - B(x,x)$ 
has a unique solution  $\overline{x} \in E$ satisfying 
$\|\overline{x}\| \leq \frac{1}{2\eta}$.
\end{Th}
\vskip 0.5cm
 \section{Proof of Theorem \ref{th1}}  
\vskip 0.5cm
In this section we shall give the proof of Theorem \ref{th1}. \\
We now need three more lemmas. In order to proceed, we define an auxiliary space $\mathcal{G}^{ \tilde q}_{q,T}$ 
which is made up of the functions $u(t,x)$ such that 
$$
\big\|u\big\|_{\mathcal{G}^{ \tilde q}_{q,T}}:= 
\Big\|\underset{0 < t < T}{\rm sup}t^{\frac{\alpha}{2}}
\big|u(t,x)\big|\Big\|_{L^{\tilde q}} < \infty,
$$
and
\begin{equation}\label{eq7}
\underset{t \rightarrow 0}{\rm lim }
\Big\|\underset{0 < \tau < t}{\rm sup}\tau^{\frac{\alpha}{2}}
\big|u(\tau,x)\big|\Big\|_{L^{\tilde q}} = 0,
\end{equation}
with 
$$
\tilde q \geq q \geq d\ {\rm and}\ \alpha 
= d\Big(\frac{1}{q} - \frac{1}{\tilde q}\Big).
$$
We recall the definition of the auxiliary space 
$\mathcal{H}^s_{p,T}$ introduced by Cannone and Planchon
\cite{M. Cannone 1999}. This space is 
made up of the functions $u(t,x)$ such that 
$$
\big\|u\big\|_{\mathcal{H}^s_{p,T}}:=
\Big\|\underset{0 < t < T}{\rm sup}
\big|\dot{\Lambda}^su(t,x)\big|\Big\|_{L^p} < \infty,
$$
and
\begin{equation}\label{eq8}
\underset{t \rightarrow 0}{\rm lim\ }
\Big\|\underset{0 < \tau < t}{\rm sup}
\big|\dot{\Lambda}^su(\tau,x)\big|\Big\|_{L^p} = 0,
\end{equation}
with
$$
p > 1\ {\rm and}\ s \geq  \frac{d}{p} - 1.
$$
The space $\mathcal{H}^s_{p,T}$ is continuously embedded 
into $L^\infty([0, T];\dot{H}^s_p(\mathbb{R}^d))$  because of the following 
elementary  inequality
$$
\underset{0 < t < T}{\rm sup}
\big\|\dot{\Lambda}^su(t,x)\big\|_{L^p} \leq 
\Big\|\underset{0 < t < T}{\rm sup}
\big|\dot{\Lambda}^su(t,x)\big|\Big\|_{L^p}.
$$
\begin{Lem}\label{lem7} Suppose that 
$u_0 \in \dot{H}^s_p(\mathbb{R}^d)$ with $p >1$ 
and $\frac{d}{p} - 1 \leq s < \frac{d}{p}$.
Then for all $\tilde q $ satisfying
$$
\tilde q > {\rm max}\{p,q\},
$$
where
$$
\frac{1}{q} = \frac{1}{p} - \frac{s}{d},
$$
we have
$$
\me^{t\Delta}u_0 \in \mathcal{G}^{\tilde q }_{q,\infty}.
$$
\end{Lem}
\begin{proof} First, we consider the case $p \leq q$.
 In this case $s \geq 0$, applying Lemma \ref{lem3} to obtain $u_0 \in L^q$. We will prove that
$$
\Big\|\underset{0 < t < \infty}{\rm sup}
t^{\frac{\alpha}{2}}\big|\me^{t\Delta}u_0\big|\Big\|_{L^{\tilde q}}
 \lesssim \big\|u_0\big\|_{L^q},\ {\rm for\ all}\ \tilde q > q.
$$
Indeed, we have the following estimates
\begin{gather}
t^{\frac{\alpha}{2}}\big|\me^{t\Delta}u_0\big| = 
\Big|\frac{t^{\frac{\alpha}{2}}}{(4\pi t)^{d/2}}
\me^{\frac{-| . |^2}{4t}}*u_0\Big| \lesssim 
\frac{1}{\sqrt{t}^{(d-\alpha)}}\me^{\frac{-| . |^2}{4t}}*|u_0| \notag \\
=  \frac{1}{|.|^{d-\alpha}}\Big(\frac{|.|}{\sqrt{t}}\Big)^{d -\alpha}
\me^{\frac{-| . |^2}{4t}}*|u_0| \leq
\underset{x \in \mathbb R^d}{\rm sup}\big(|x|^{d-\alpha}
\me^{\frac{-|x|^2}{4}}\big).\frac{1}{|.|^{d-\alpha}}*\big|u_0\big| \notag \\ 
\lesssim \frac{1}{|.|^{d-\alpha}}*\big|u_0\big|. \label{eq9}  
\end{gather}
From the estimate \eqref{eq9}, applying Lemma \ref{lem5}(c) to obtain
\begin{gather*}
\Big\|\underset{0 < t < \infty}{\rm sup}t^{\frac{\alpha}{2}}
\big|\me^{t\Delta}u_0\big|\Big\|_{L^{\tilde q}} 
\lesssim \Big\|\frac{1}{|.|^{d-\alpha}}*\big|u_0\big|\Big\|_{L^{\tilde q}}
 \lesssim \Big\|\frac{1}{|.|^{d-\alpha}}\Big\|_{L^{\frac{d}{d-\alpha}, \infty}}
\big\|u_0\big\|_{L^{q,\tilde q}}\notag  \\  
\lesssim \big\|u_0\big\|_{L^q},\ ({\rm note\ that}\ \frac{1}{|.|^s}
\ \in L^{\frac{d}{s},\infty}(\mathbb{R}^d)\ {\rm with}\ 0 < s \leq d).
\end{gather*}
This proves the result. We now prove that
$$
\underset{t \rightarrow 0}{\rm lim }
\Big\|\underset{0 < \tau < t}{\rm sup}
\tau^{\frac{\alpha}{2}}\big|\me^{\tau\Delta}
u_0\big|\Big\|_{L^{\tilde q}} = 0.
$$
Set $\mathcal{X}_n(x) = 0$ for $x \in \{x : \ |x| < n\}
 \cap \{x : \ |u_0(x)| < n\}$ and $\mathcal{X}_n(x) = 1$ otherwise. 
We have
\begin{equation}\label{eq10}
t^{\frac{\alpha}{2}}\big|\me^{t\Delta}u_0\big| \leq  
C\Big(\frac{1}{\sqrt{t}^{(d-\alpha)}}
\me^{\frac{-| . |^2}{4t}}*|\mathcal{X}_nu_0| 
+ \frac{1}{\sqrt{t}^{(d-\alpha)}}
\me^{\frac{-| . |^2}{4t}}*|(1 - \mathcal{X}_n)u_0| \Big).
\end{equation}
Let $\hat q$ be fixed such that $q < \hat q < \tilde q$ 
and $\beta = d(\frac{1}{\hat q} - \frac{1}{\tilde q})$. 
Arguing as in the proof of the estimate \eqref{eq9}, we derive
\begin{gather}
C\frac{1}{\sqrt{t}^{(d-\alpha)}}
\me^{\frac{-| . |^2}{4t}}*|\mathcal{X}_nu_0| 
\leq C_1\frac{1}{|.|^{d-\alpha}}*
\big|\mathcal{X}_nu_0\big|, \label{eq11}
\end{gather}
and
\begin{gather}
C\frac{1}{\sqrt{t}^{(d-\alpha)}}
\me^{\frac{-| . |^2}{4t}}*|(1 - \mathcal{X}_n)u_0| =  
Ct^{\frac{\alpha -\beta}{2}}\frac{1}{\sqrt{t}^{(d-\beta)}}
\me^{\frac{-| . |^2}{4t}}*|(1 - \mathcal{X}_n)u_0| \notag \\
 \leq C\underset{x \in \mathbb R^d}{\rm sup}(|x|^{d-\beta}
\me^{\frac{-|x|^2}{4}})t^{\frac{\alpha -\beta}{2}}
\frac{1}{|.|^{d-\beta}}*\big|(1 - \mathcal{X}_n)u_0\big| \notag \\ \leq 
C_2nt^{\frac{d}{2}(\frac{1}{q} - \frac{1}{\hat q})}
\frac{1}{|.|^{d-\beta}}*\big|1 - \mathcal{X}_n\big|. \label{eq12}
\end{gather}
From the estimates \eqref{eq10}, \eqref{eq11}, 
and \eqref{eq12}, we have
\begin{gather}
\Big\|\underset{0 < \tau < t}{\rm sup}\tau^{\frac{\alpha}{2}}
\big|\me^{\tau\Delta}u_0\big|\Big\|_{L^{\tilde q}} \leq  \notag \\
C_1\Big\|\frac{1}{|.|^{d-\alpha}}*\big|\mathcal{X}_nu_0
\big|\Big\|_{L^{\tilde q}} + C_2nt^{\frac{d}{2}(\frac{1}{q} 
- \frac{1}{\hat q})}\Big\|\frac{1}{|.|^{d-\beta}}*\big|1 - \mathcal{X}_n
\big|\Big\|_{L^{\tilde q}}\leq  \notag\\
C_3\Big\|\frac{1}{|.|^{d-\alpha}}\Big\|_{L^{\frac{d}{d-\alpha}, \infty}}
\big\|\mathcal{X}_nu_0\big\|_{L^q} + C_4nt^{\frac{d}{2}(\frac{1}{q} 
- \frac{1}{\hat q})}\Big\|\frac{1}{|.|^{d-\beta}}
\Big\|_{L^{\frac{d}{d-\beta}, \infty}}\big\|1 
- \mathcal{X}_n\big\|_{L^{\hat q}}\leq \notag \\
 C_5\big\|\mathcal{X}_nu_0\big\|_{L^q} 
+ C_6nt^{\frac{d}{2}(\frac{1}{q} - \frac{1}{\hat q})}
\big\|1 - \mathcal{X}_n\big\|_{L^{\hat q}}. \label{eq13}
\end{gather}
For any $\epsilon > 0$, we can take $n$ large enough that 
\begin{equation}\label{eq14}
C_5\big\|\mathcal{X}_nu_0\big\|_{L^q} < \frac{\epsilon}{2}.
\end{equation}
Fixed one of such $n$, there exists $t_0=t_0(n)> 0$ satisfying 
\begin{equation}\label{eq15}
C_6nt^{\frac{d}{2}(\frac{1}{q} - \frac{1}{\hat q})}
\big\|1 - \mathcal{X}_n\big\|_{L^{\hat q}} < \frac{\epsilon}{2},\ {\rm for}\  t<t_0.
\end{equation} 
From the estimates 
\eqref{eq13}, \eqref{eq14}, and \eqref{eq15}, we have
$$
\Big\|\underset{0 < \tau < t}{\rm sup}\tau^{\frac{\alpha}{2}}
\big|\me^{\tau\Delta}u_0\big|\Big\|_{L^{\tilde q}} \leq C_5
\big\|\mathcal{X}_nu_0\big\|_{L^q} + 
C_6nt^{\frac{d}{2}(\frac{1}{q} - \frac{1}{\hat q})}\big\|1 
- \mathcal{X}_n\big\|_{L^{\hat q}} < \epsilon, {\rm for}\ t<t_0.
$$
We now consider the case $p > q$. In this case $s <0$. We prove that
$$
\Big\|\underset{0 < t < \infty}{\rm sup}t^{\frac{\alpha}{2}}
\big|\me^{t\Delta}u_0\big|\Big\|_{L^{\tilde q}} 
\lesssim \big\|u_0\big\|_{\dot{H}^s_p},\ {\rm for\  all}\ \tilde q > p.
$$
We have 
\begin{equation*}
\me^{t\Delta}u_0 =\me^{t\Delta}\dot{\Lambda}^{-s}\dot{\Lambda}^su_0 
= \frac{1}{t^{\frac{d - s}{2}}} K\Big(\frac{.}{\sqrt t}\Big) * 
(\dot{\Lambda}^su_0),
\end{equation*}
where
$$
\hat K(\xi) = \frac{1}{(2\pi)^{\frac{d}{2}}} \me^{-|\xi|^2} |\xi|^{- s}\ {\rm and}\  |K(x)| \lesssim \frac{1}{(1 + |x|)^{d - s}}.
$$
From the above inequality, we have
\begin{gather}
t^{\frac{\alpha}{2}}\big|\me^{t\Delta}u_0\big| 
\leq \Big|\frac{1}{\sqrt{t}^{(d-s-\alpha)}}K\Big(\frac{.}{\sqrt t}\Big)
\Big|*|\dot{\Lambda}^su_0| \notag \\ 
=  \Big|\frac{1}{|.|^{d-s-\alpha}}\Big(\frac{|.|}{\sqrt t}\Big)^{d-s-\alpha}
K\Big(\frac{.}{\sqrt t}\Big)\Big|*|\dot{\Lambda}^su_0| \notag \\
\leq \underset{x \in \mathbb R^d}{\rm sup}\big(\big| |x|^{d-s-\alpha}
K(x)\big|\big)\frac{1}{|.|^{d-s-\alpha}}*\big|\dot{\Lambda}^su_0\big| 
\lesssim \frac{1}{|.|^{d-s-\alpha}}*
\big|\dot{\Lambda}^su_0\big|. \label{eq16}  
\end{gather}
From the estimate \eqref{eq16}, applying Lemma \ref{lem5}(c), we have 
\begin{gather*}
\Big\|\underset{0 < t < \infty}{\rm sup}t^{\frac{\alpha}{2}}
\big|\me^{t\Delta}u_0\big|\Big\|_{L^{\tilde q}} \lesssim 
\Big\|\frac{1}{|.|^{d-s-\alpha}}*\big|\dot{\Lambda}^su_0
\big|\Big\|_{L^{\tilde q}} \leq  \\  
\Big\|\frac{1}{|.|^{d-s-\alpha}}*
\big|\dot{\Lambda}^su_0\big|\Big\|_{L^{\tilde q,p}} 
\lesssim \Big\|\frac{1}{|.|^{d-s-\alpha}}
\Big\|_{L^{\frac{d}{d-s-\alpha}, \infty}}\big\|\dot{\Lambda}^s
u_0\big\|_{L^p} \simeq \big\|u_0\big\|_{\dot{H}^s_p}.
\end{gather*}
This proves the result. We now claim that
\begin{equation*}
\underset{t \rightarrow 0}{\rm lim\ }
\Big\|\underset{0 < \tau < t}{\rm sup}\tau^{\frac{\alpha}{2}}
\big|\me^{\tau\Delta}u_0\big|\Big\|_{L^{\tilde q}} = 0,
\ \text{ for all}\ \tilde q > p.
\end{equation*}
Set  $\mathcal{X}_{n,s}(x) = 0$ for $x \in \{x : \ |x| < n\} 
\cap \{x : \ |\dot{\Lambda}^su_0(x)| < n\}$ and 
$\mathcal{X}_{n,s}(x) = 1$ otherwise. Let $\hat q$ 
be fixed such that $p < \hat q <\tilde q$ and $\beta = d(\frac{1}{p} 
- \frac{1}{\hat q})$. For any $\epsilon > 0$, by an arguing similar to the case $q > p$, there exist a sufficiently large $n$ and a sufficiently small $t_0 = t_0(n)$ such that
\begin{gather*}
\Big\|\underset{0 < \tau < t}{\rm sup}\tau^{\frac{\alpha}{2}}
\big|\me^{\tau\Delta}u_0\big|\Big\|_{L^{\tilde q}} \leq 
C_1\Big\|\frac{1}{|.|^{d-s-\alpha}}\Big\|_{L^{\frac{d}{d-s
-\alpha}, \infty}}\big\|\mathcal{X}_{n,s}\dot{\Lambda}^s
u_0\big\|_{L^p} \notag \\
 + \ C_2nt^{\frac{\beta}{2}}\Big\|\frac{1}{|.|^{d-s-\alpha+\beta}}
\Big\|_{L^{\frac{d}{d-s-\alpha+\beta}, \infty}}\big\|1
 - \mathcal{X}_{n,s}\big\|_{L^{\hat q}} < \epsilon,\ {\rm for}\ t<t_0.
\end{gather*}
\end{proof}
In the following lemmas a particular attention will be 
devoted to the study of the bilinear operator $B(u, v)(t)$ defined by 
\begin{equation}\label{eq17}
B(u, v)(t) = \int_{0}^{t} e^{(t-\tau ) \Delta} \mathbb{P} 
\nabla.\big(u(\tau)\otimes v(\tau)\big) \dif\tau.
\end{equation}
\begin{Lem}\label{lem8}
Let $p$ and $s$ be such that 
$$
p > \frac{d}{2}\ \text{and}\ \frac{d}{p} - 1 \leq s  < \frac{d}{2p}.
$$
Then the bilinear operator $B$ is continuous 
from $\mathcal{G}^{\tilde q}_{q,T} \times 
\mathcal{G}^{\tilde q}_{q,T}$ into 
$\mathcal{H}^s_{p,T}$, where 
$$
\frac{1}{q} = \frac{1}{p} - \frac{s}{d},\  q < \tilde q  <  2p,
$$
and we have the inequality
\begin{equation}\label{eq18}
\big\|B(u, v)\big\|_{\mathcal{H}^s_{p,T}} \leq
 CT^{\frac{1}{2}(1 + s  - \frac{d}{p})}
\big\|u\big\|_{\mathcal{G}^{\tilde q}_{q,T}}
\big\|v\big\|_{\mathcal{G}^{\tilde q}_{q,T}},
\end{equation}
where C is a positive constant and independent of T.
\end{Lem}
\begin{proof}From the equality \eqref{eq17}, applying Lemma \ref{lem2} to obtain
\begin{gather}
\big|\dot{\Lambda}^sB(u,v)(t)(x)\big| \leq \int_{0}^{t}
 \big|\dot{\Lambda}^s\me^{(t- \tau) \Delta} \mathbb{P} \nabla .
\big(u(\tau,x)\otimes v(\tau,x)\big)\big|{\rm d}\tau  \notag \\
= \int_{0}^{t} \Big|K_{t - \tau}(x)*\big(u(\tau,x)
 \otimes v(\tau,x)\big)\Big|{\rm d}\tau \notag \\ 
\lesssim  \int_{0}^{t} \Big|\frac{1}{{(t - \tau)}^{\gamma_2}.|x|^{\gamma_1}}*\big(u(\tau,x)
 \otimes v(\tau,x)\big)\Big|{\rm d}\tau \label{eq19}
\end{gather}
where
$$
\gamma_1 >0, \gamma_2 > 0, \gamma_1 + 2\gamma_2 = d+1+s.
$$
Using the estimate \eqref{eq19} for
$$
\gamma_1 = d\Big(1 +\frac{1}{p} - \frac{2}{\tilde q}\Big), 
\gamma_2 = \frac{1}{2} -\frac{d}{2p} 
+ \frac{s}{2} + \frac{d}{\tilde q}, 
$$
and applying Lemma \ref{lem6} to obtain
\begin{gather}
\big|\dot{\Lambda}^sB(u,v)(t)(x)\big| \notag \\ 
\lesssim  \frac{1}{|x|^{d(1 +\frac{1}{p} - \frac{2}{\tilde q})}}*
\int_0^t\frac{1}{(t - \tau)^{\frac{1}{2} -\frac{d}{2p} 
+ \frac{s}{2} + \frac{d}{\tilde q}}}\big|u(\tau,x) 
\otimes v(\tau,x)\big|{\rm d}\tau \notag \\ 
\lesssim  \frac{1}{|x|^{d(1 +\frac{1}{p} - \frac{2}{\tilde q})}}*
\int_0^t\frac{1}{(t - \tau)^{\frac{1}{2} -\frac{d}{2p} 
+ \frac{s}{2} + \frac{d}{\tilde q}}}\tau^{-\alpha}
\underset{0 < \eta < t}{\rm sup}\eta^{\frac{\alpha}{2}}
\big|u(\eta,x)\big|\underset{0 < \eta < t}{\rm sup}
\eta^{\frac{\alpha}{2}}\big|v(\eta,x)\big|{\rm d}\tau \notag\\
= \frac{1}{|x|^{d(1 +\frac{1}{p} - \frac{2}{\tilde q})}}*
\Big(\underset{0 < \eta < t}{\rm sup}\eta^{\frac{\alpha}{2}}
\big|u(\eta,x)\big|\underset{0 < \eta < t}{\rm sup}
\eta^{\frac{\alpha}{2}}\big|v(\eta,x)\big|\Big)\int_0^t
\frac{1}{(t - \tau)^{\frac{1}{2} -\frac{d}{2p} 
+ \frac{s}{2} + \frac{d}{\tilde q}}}\tau^{-\alpha}{\rm d}\tau  \notag \\ 
\simeq t^{\frac{1}{2}(1 + s  - \frac{d}{p})}\frac{1}{|x|^{d(1 +\frac{1}{p}
 - \frac{2}{\tilde q})}}*\Big(\underset{0 < \eta < t}{\rm sup}
\eta^{\frac{\alpha}{2}}\big|u(\eta,x)\big|\underset{0 < \eta < t}{\rm sup}
\eta^{\frac{\alpha}{2}}\big|v(\eta,x)\big|\Big).
  \label{eq20}
\end{gather}
From the estimate \eqref{eq20}, applying Lemma \ref{lem5}(c) and H\"{o}lder's inequality 
in order to obtain 
\begin{gather}
\Big\|\underset{0 < \tau < t}{\rm sup}\big|\dot{\Lambda}^s
B(u,v)(\tau)\big|\Big\|_{L^p} \leq \Big\|\underset{0 < \tau < t}
{\rm sup}\big|\dot{\Lambda}^sB(u,v)(\tau)\big|\Big\|_{L^{p, 
\frac{\tilde q}{2}}} \notag \\
\lesssim t^{\frac{1}{2}(1 + s  - \frac{d}{p})}
\Big\|\frac{1}{|x|^{d(1 +\frac{1}{p} - \frac{2}{\tilde q})}}
\Big\|_{L^{\frac{1}{1+ \frac{1}{p} - \frac{2}{\tilde q}},\infty}}
\Big\|\underset{0 < \eta < t}{\rm sup}\eta^{\frac{\alpha}{2}}
\big|u(\eta,x)\big|\underset{0 < \eta < t}{\rm sup}
\eta^{\frac{\alpha}{2}}\big|v(\eta,x)\big|\Big\|_{L^{\frac{\tilde q}{2}}} \notag \\  
\lesssim t^{\frac{1}{2}(1 + s  - \frac{d}{p})}
\Big\|\underset{0 < \eta < t}{\rm sup}\eta^{\frac{\alpha}{2}}
\big|u(\eta,x)\big|\Big\|_{L^{\tilde q}}
\Big\|\underset{0 < \eta < t}{\rm sup}
\eta^{\frac{\alpha}{2}}\big|v(\eta,x)
\big|\Big\|_{L^{\tilde q}}. \label{eq21}
\end{gather}
Let us now check the validity of the condition \eqref{eq8} for the bilinear term $B(u,v)(t)$.
In fact, from the estimate \eqref{eq21} it follows that
\begin{equation*}
\underset{t \rightarrow 0}{\rm lim }
\Big\|\underset{0 < \tau < t}{\rm sup}
\big|\dot{\Lambda}^sB(u,v)(\tau)\big|\Big\|_{L^p} = 0,
\end{equation*}
whenever
$$
\underset{t \rightarrow 0}{\rm lim\ }
\Big\|\underset{0 < \tau < t}{\rm sup}\tau^{\frac{\alpha}{2}}
\big|u(\tau,x)\big|\Big\|_{L^{\tilde q}} 
= \underset{t \rightarrow 0}{\rm lim }
\Big\|\underset{0 < \tau < t}{\rm sup}
\tau^{\frac{\alpha}{2}}\big|v(\tau,x)\big|\Big\|_{L^{\tilde q}} = 0.
$$
The estimate \eqref{eq18} is deduced from the inequality \eqref{eq21}.\end{proof} 
\begin{Lem}\label{lem9} Let $q$ and $q_1$ be such that $d\leq q < q_1 <+\infty$. Then the bilinear operator $B$ is continuous from  $\mathcal{G}^{q_1}_{q,T} \times \mathcal{G}^{q_1}_{q,T}$ into $\mathcal{G}^{q_2}_{q,T}$ for all $q_2$ satisfying $
\frac{1}{q_2} \in \big(0,\frac{1}{q}\big] \cap \big(\frac{2}{q_1}-\frac{1}{d},\ \frac{2}{q_1}\big)$, and we have the inequality
\begin{equation}\label{eq22}
\big\|B(u, v)\big\|_{\mathcal{G}^{q_2}_{q,T}} \leq
 CT^{\frac{1}{2}(1  - \frac{d}{q})}
\big\|u\big\|_{\mathcal{G}^{q_1}_{q,T}}
\big\|v\big\|_{\mathcal{G}^{q_1}_{q,T}},
\end{equation}
where C is a positive constant and independent of T.
\end{Lem}
\begin{proof} From the equality \eqref{eq17}, applying Lemma \ref{lem2} to obtain
\begin{gather}
\big|B(u,v)(t)(x)\big| \leq \int_{0}^{t}
 \big|\me^{(t- \tau) \Delta} \mathbb{P} \nabla .
\big(u(\tau,x)\otimes v(\tau,x)\big)\big|{\rm d}\tau  \notag \\
= \int_{0}^{t} \Big|K_{t - \tau}(x)*\big(u(\tau,x)
 \otimes v(\tau,x)\big)\Big|{\rm d}\tau \notag \\ 
\lesssim  \int_{0}^{t} \Big|\frac{1}{{(t - \tau)}^{\gamma_2}.|x|^{\gamma_1}}*\big(u(\tau,x)
 \otimes v(\tau,x)\big)\Big|{\rm d}\tau \label{eq23}
\end{gather}
where
$$
\gamma_1 >0, \gamma_2 > 0, \gamma_1 + 2\gamma_2 = d+1.
$$
Set
$$
\alpha_1 =d(\frac{1}{q}-\frac{1}{q_1}), \alpha_2 =d(\frac{1}{q}-\frac{1}{q_2}).
$$
Using the estimate \eqref{eq23} for
$$
\gamma_1 = d\Big(1 +\frac{1}{q_2} - \frac{2}{q_1}\Big), 
\gamma_2 = \frac{1}{2} -\frac{d}{2q_2} + \frac{d}{q_1}, 
$$
and applying Lemma \ref{lem6} to obtain
\begin{gather}
\big|B(u,v)(t)(x)\big|\lesssim  \frac{1}{|x|^{d(1 +\frac{1}{q_2} - \frac{2}{q_1})}}*
\int_0^t\frac{1}{(t - \tau)^{\frac{1}{2} -\frac{d}{2q_2} + \frac{d}{q_1}}}\big|u(\tau,x) 
\otimes v(\tau,x)\big|{\rm d}\tau \notag \\ 
\lesssim  \frac{1}{|x|^{d(1 +\frac{1}{q_2} - \frac{2}{q_1})}}*
\int_0^t\frac{1}{(t - \tau)^{\frac{1}{2} -\frac{d}{2q_2} + \frac{d}{q_1}}}\tau^{-\alpha_1}
\underset{0 < \eta < t}{\rm sup}\eta^{\frac{\alpha_1}{2}}
\big|u(\eta,x)\big|\underset{0 < \eta < t}{\rm sup}
\eta^{\frac{\alpha_1}{2}}\big|v(\eta,x)\big|{\rm d}\tau \notag\\
= \frac{1}{|x|^{d(1 +\frac{1}{q_2} - \frac{2}{q_1})}}*
\Big(\underset{0 < \eta < t}{\rm sup}\eta^{\frac{\alpha_1}{2}}
\big|u(\eta,x)\big|\underset{0 < \eta < t}{\rm sup}
\eta^{\frac{\alpha_1}{2}}\big|v(\eta,x)\big|\Big)\int_0^t
\frac{1}{(t - \tau)^{\frac{1}{2} -\frac{d}{2q_2} + \frac{d}{q_1}}\tau^{\alpha_1}}{\rm d}\tau  \notag \\ 
\simeq t^{\frac{1}{2}(1 - \frac{d}{q})-\frac{\alpha_2}{2}}\frac{1}{|x|^{d(1 +\frac{1}{q_2} - \frac{2}{q_1})}}*\Big(\underset{0 < \eta < t}{\rm sup}
\eta^{\frac{\alpha_1}{2}}\big|u(\eta,x)\big|\underset{0 < \eta < t}{\rm sup}
\eta^{\frac{\alpha_1}{2}}\big|v(\eta,x)\big|\Big).
  \label{eq24}
\end{gather}
From the estimate \eqref{eq24}, applying Lemma \ref{lem5}(c) and H\"{o}lder's inequality 
in order to obtain 
\begin{gather}
\Big\|\underset{0 < \tau < t}{\rm sup}t^{\frac{\alpha_2}{2}}\big|
B(u,v)(\tau)\big|\Big\|_{L^{q_2}} \leq \Big\|\underset{0 < \tau < t}
{\rm sup}t^{\frac{\alpha_2}{2}}\big|B(u,v)(\tau)\big|\Big\|_{L^{q_2, 
\frac{q_2}{2}}} \notag \\
\lesssim t^{\frac{1}{2}(1-\frac{d}{q})}
\Big\|\frac{1}{|x|^{d(1 +\frac{1}{q_2} - \frac{2}{q_1})}}
\Big\|_{L^{\frac{1}{1 +\frac{1}{q_2} - \frac{2}{q_1}},\infty}}
\Big\|\underset{0 < \eta < t}{\rm sup}\eta^{\frac{\alpha_1}{2}}
\big|u(\eta,x)\big|\underset{0 < \eta < t}{\rm sup}
\eta^{\frac{\alpha_1}{2}}\big|v(\eta,x)\big|\Big\|_{L^{\frac{q_1}{2}}} \notag \\  
\lesssim t^{\frac{1}{2}(1-\frac{d}{q})}
\Big\|\underset{0 < \eta < t}{\rm sup}\eta^{\frac{\alpha_1}{2}}
\big|u(\eta,x)\big|\Big\|_{L^{\tilde q}}
\Big\|\underset{0 < \eta < t}{\rm sup}
\eta^{\frac{\alpha_1}{2}}\big|v(\eta,x)
\big|\Big\|_{L^{q_1}}. \label{eq25}
\end{gather}
Let us now check the validity of the condition \eqref{eq7} for the bilinear term $B(u,v)(t)$.
In fact, from the estimate \eqref{eq25} it follows that
\begin{equation*}
\underset{t \rightarrow 0}{\rm lim }
\Big\|\underset{0 < \tau < t}{\rm sup}t^{\frac{\alpha_2}{2}}\big|B(u,v)(\tau)\big|\Big\|_{L^{q_2}} = 0,
\end{equation*}
whenever
$$
\underset{t \rightarrow 0}{\rm lim\ }
\Big\|\underset{0 < \tau < t}{\rm sup}\tau^{\frac{\alpha_1}{2}}
\big|u(\tau,x)\big|\Big\|_{L^{q_1}} 
= \underset{t \rightarrow 0}{\rm lim }
\Big\|\underset{0 < \tau < t}{\rm sup}
\tau^{\frac{\alpha_1}{2}}\big|v(\tau,x)\big|\Big\|_{L^{q_1}} = 0
$$
The estimate \eqref{eq22} is deduced from the inequality \eqref{eq25}.\end{proof} 

\vskip 0.5cm
{\bf Proof of Theorem \ref{th1}}
\vskip 0.5cm
(a) Applying Lemma \ref{lem9} for $q_1=q_2 = \tilde q$, we deduce that $B$ is 
continuous from $\mathcal{G}^{\tilde q}_{q,T} 
\times \mathcal{G}^{\tilde q}_{q,T}$ to 
$\mathcal{G}^{\tilde q}_{q,T}$ and we have the inequality
$$
\big\|B(u, v)\big\|_{\mathcal{G}^{\tilde q}_{q,T}} 
\leq C_{q,\tilde q, d}T^{\frac{1}{2}(1-\frac{d}{q})}
\big\|u\big\|_{\mathcal{G}^{\tilde q}_{q,T}}
\big\|v\big\|_{\mathcal{G}^{\tilde q}_{q,T}} =
 C_{q,\tilde q, d}T^{\frac{1}{2}(1+s-\frac{d}{p})}
\big\|u\big\|_{\mathcal{G}^{\tilde q}_{q,T}}
\big\|v\big\|_{\mathcal{G}^{\tilde q}_{q,T}},
$$
where $C_{q,\tilde q,d}$ is a positive  and independent 
of  $T$. From Theorem \ref{th2} and the above inequality, 
we deduce that for any $u_0 \in  \dot{H}^s_p$ satisfying
$$
T^{\frac{1}{2}(1+s-\frac{d}{p})}
\big\|\me^{t\Delta}u_0\big\|_{\mathcal{G}^{\tilde q}_{q,T}} 
= T^{\frac{1}{2}(1+s-\frac{d}{p})}\underset{0 < t < T}{{\rm sup}}
t^{\frac{\alpha}{2}}\big\|\me^{t\Delta}u_0\big\|_{L^{\tilde q}} 
\leq  \frac{1}{4C_{q,\tilde q, d}},
$$
where
$$
\alpha = d\Big(\frac{1}{q}- \frac{1}{\tilde q}\Big) 
= d\Big(\frac{1}{p}- \frac{s}{d} - \frac{1}{\tilde q}\Big),
$$
NSE has a solution $u$ on the interval $(0, T)$ 
so that $u \in \mathcal{G}^{\tilde q}_{q,T}$. We prove that  $u \in \underset{r > {\rm max}\{p, q\}}
{\bigcap}\mathcal{G}^r_{q,T}$.  Indeed, applying Lemma \ref{lem9}, we have $B(u,u) \in  \mathcal{G}^{r}_{q,T}$ for all $r$ satisfying
$\frac{1}{r} \in \big(0,\frac{1}{q}\big] \cap \big(\frac{2}{\tilde q}-\frac{1}{d},\frac{2}{\tilde q}\big)$. Applying Lemma \ref{lem7}, we have
$e^{t\Delta}u_0 \in \mathcal{G}^{r}_{q,T}$ for all  $r$ satisfying $\frac{1}{r} \in \big(0, \frac{1}{{\rm max}\{p,q\}}\big)$. Since $u = e^{t\Delta}u_0 -B(u,u)$, it follows that $u \in \mathcal{G}^{r}_{q,T}$ for all  $r$ satisfying  $\frac{1}{r} \in \big(0, \frac{1}{{\rm max}\{p,q\}}\big)\cap\big(\frac{2}{\tilde q}-\frac{1}{d},\frac{2}{\tilde q}\big)$. Applying again Lemmas \ref{lem9} and \ref{lem7}, in exactly the same way, since $u \in \mathcal{G}^{r}_{q,T}$ for all $r$ satisfying  $\frac{1}{r} \in \big(0, \frac{1}{{\rm max}\{p,q\}}\big)\cap\big(\frac{2}{\tilde q}-\frac{1}{d},\frac{2}{\tilde q}\big)$,
 it follows that $u \in \mathcal{G}^{r}_{q,T}$ for all $r$ satisfying $\frac{1}{r} \in \big(0, \frac{1}{{\rm max}\{p,q\}}\big)\cap\big(\frac{1}{d}-2^2(\frac{1}{d}-\frac{1}{\tilde q}), \frac{2^2}{\tilde q}\big)$. By induction, 
we get $u \in \mathcal{G}^{r}_{q,T}$ for all $r$ satisfying 
$\frac{1}{r} \in \big(0, \frac{1}{{\rm max}\{p,q\}}\big)\cap\big(\frac{1}{d}-2^n(\frac{1}{d}-\frac{1}{\tilde q}),\frac{2^n}{\tilde q}\big)$ with $n\geq 1$. Since $\frac{1}{d}-\frac{1}{\tilde q} >0$, it follows that there exists sufficiently large $n$ satisfying $ \big(0, \frac{1}{{\rm max}\{p,q\}}\big)\cap\big(\frac{1}{d}-2^n(\frac{1}{d}-\frac{1}{\tilde q}),\frac{2^n}{\tilde q}\big) = \big(0, \frac{1}{{\rm max}\{p,q\}}\big)$ . Therefore
$u \in \mathcal{G}^r_{q,T}$ for all $r > {\rm max}\{p,q\}$. This proves the result.\\
We now prove that $u \in L^\infty([0, T]; \dot{H}^s_p)$. 
Indeed, from $u \in \mathcal{G}^r_{q,T}$ for all $r > {\rm max}\{p,q\}$, applying Lemma \ref{lem8} to obtain $B(u,u) \in \mathcal{H}^s_{p,T} \subseteq L^\infty\big([0, T]; \dot{H}^s_p\big)$. On the other hand, since $u \in \dot{H}^s_p$, it follows that  $\me^{t\Delta}u_0 \in L^\infty\big([0, T]; \dot{H}^s_p\big)$. Therefore
$$
u = \me^{t\Delta}u_0 - B(u,u) \in  L^\infty\big([0, T]; \dot{H}^s_p\big).
$$
Finally, we will show that the condition \eqref{eq1} is valid when $T$ is small enough. From the definition of $\mathcal{G}^{ \tilde q}_{q,T}$ and Lemma \ref{lem7}, we deduce that the left-hand side of the condition \eqref{eq1} converges to $0$ when $T$ goes to $0$. Therefore the condition \eqref{eq1} holds for arbitrary $u_0 \in \dot{H}^s_p (\mathbb{R}^d)$ when $T(u_0)$ is small enough. \\
(b) From Lemma \ref{lem4}, the two quantities $\big\|u_0\big\|_{\dot{F}^{\frac{d}{\tilde q} - 1, \infty}_{\tilde q}}$ and  $\underset{0 < t < \infty}{{\rm sup}}t^{\frac{1}{2}(1- \frac{d}{\tilde q})}\big\|\me^{t\Delta}u_0\big\|_{L^{\tilde q}}$  are equivalent. Thus, there exists a positive constant $\sigma_{\tilde q,d}$ such that  the condition \eqref{eq1} holds for  $T = \infty$ whenever $\big\|u_0\big\|_{\dot{F}^{\frac{d}{\tilde q} - 1, \infty}_{\tilde q}} \leq \sigma_{\tilde q,d}$.\qed 
\vskip 0.5cm
{\bf Proof of Proposition \ref{pro1}}
\vskip 0.5cm
By Theorem \ref{th1}, we only need to prove that $w \in \mathcal{H}^{\frac{d}{\tilde p}-1}_{{\tilde p},T}$ for all $\tilde p >  \frac{1}{2}{\rm max}\{p, d\}$ and if $p \geq d$ then $\underset{0 < t < T}{\rm sup}t^{\frac{1}{2}({1-\frac{d}{p}})}\big|u(t, .)\big| \in L^p$. Indeed, applying Lemma \ref{lem8}, we deduce that the bilinear operator $B$ is continuous from $\mathcal{G}^r_{d,T} \times \mathcal{G}^r_{d,T}$ into $\mathcal{H}^{\frac{d}{\tilde p}-1}_{{\tilde p},T}$ for all $\tilde p > \frac{d}{2}$ and $r$ satisfying  $d < r < 2\tilde p$; hence from $u \in \underset{r > {\rm max}\{p,d\}}{\bigcap}\mathcal{G}^r_{d,T}$ and $2\tilde p > {\rm max}\{p,d\}$, we have $w=-B(u,u) \in \mathcal{H}^{\frac{d}{\tilde p}-1}_{{\tilde p},T}$.\\
We now prove that if $p \geq d$ then $\underset{0 < t < T}{\rm sup}t^{\frac{1}{2}({1-\frac{d}{p}})}\big|u(t, .)\big| \in L^p$.
Indeed, we notice that, if $u_0 \in \dot{H}^{\frac{d}{p}-1}_p$, 
then $u_0 = \dot{\Lambda}^{1-\frac{d}{p}}v_0$ 
with $v_0 \in L^p$; hence $t^{\frac{1}{2}({1-\frac{d}{p}})}
\big|\me^{t\Delta}u_0\big| \lesssim M_{v_0}$, where $M_{v_0}$ 
is the Hardy-Littlewood maximal function of $v_0$
(hence $M_{v_0} \in L^p$). On the other hand, 
from $u \in \underset{r > p}{\bigcap}\mathcal{G}^r_{d,T}$, 
we apply   Lemma \ref{lem9} to obtain $B(u,u) 
\in \mathcal{G}^p_{d,T}$, hence $\underset{0 < t < T}{\rm sup}t^{\frac{1}{2}({1-\frac{d}{p}})}\big|B(u,u)(t, .)\big| \in L^p$. 
Thus
$$
\underset{0 < t < T}{\rm sup}t^{\frac{1}{2}({1-\frac{d}{p}})}
\big|u(t, .)\big|  \leq \underset{0 < t <T}{\rm sup}
t^{\frac{1}{2}({1-\frac{d}{p}})}\big|\me^{t\Delta}u_0\big| 
+ \underset{0 < t < T}{\rm sup}t^{\frac{1}{2}({1-\frac{d}{p}})}
\big|B(u,u)(t, .)\big| \in L^p.
$$
\qed
\vskip 0.5cm
{\bf Proof of Proposition \ref{pro2}}
\vskip 0.5cm
By Lemma \ref{lem4}, we deduce that two quantities 
$$
\big\|u_0\big\|_{\dot{F}^{s- (\frac{d}{p} 
- \frac{d}{\tilde q}), \infty}_{\tilde q}}\ {\rm and}\ \big\|\underset{0 < t < \infty}{\rm sup}
t^{\frac{d}{2}(\frac{1}{p}- \frac{s}{d} 
- \frac{1}{\tilde q})}|\me^{t\Delta}u_0|\big\|_{L^{\tilde q}}
$$ 
are equivalent. Thus
$$
\Big\|\underset{0 < t < T}{\rm sup}
t^{\frac{d}{2}(\frac{1}{p}- \frac{s}{d} - \frac{1}{\tilde q})}
\big|\me^{t\Delta}u_0\big|\Big\|_{L^{\tilde q}} 
\lesssim \big\|u_0\big\|_{\dot{F}^{s- (\frac{d}{p} 
- \frac{d}{\tilde q}), \infty}_{\tilde q}}.
$$
The Proposition \ref{pro2} is proved by applying 
Theorem \ref{th1} and the above inequality. \qed

\subsection*{Acknowledgements} This research is funded by 
Vietnam National  Foundation for Science and Technology 
Development (NAFOSTED) under grant number  101.02-2014.50.

\end{document}